\documentclass[11pt]{amsart}
\numberwithin{equation}{section}
\newtheorem{theorem}{Theorem}

\numberwithin{theorem}{section}

\def\al{\aligned}
\def\eal{\endaligned}
\def\M{{\bf M}}
\def\be{\begin{equation}}
\def\ee{\end{equation}}
\def\lab{\label}

\def\M{{\bf M}}

\def\al{\aligned}

\def\p{\partial}

\def\lam{\lambda}
\numberwithin{equation}{section}

\begin{document}

\tracingpages 1
\title[nonsingular]{\bf A note on convergence of noncompact nonsingular solutions of the Ricci flow
}
\author{Qi S. Zhang}
\address{
Department of
Mathematics, University of California, Riverside, CA 92521, USA}
\date{2020.09; MSC2020: 58J35.}

\begin{abstract}
We extend some convergence results on nonsingular compact Ricci flows in the papers \cite{Ha:1}, \cite{Se:1} and \cite{FZZ:2} to certain infinite volume
noncompact cases which are "partially" nonsingular. As an application, for  a finite time singularity which is partially type I,  it is shown that a blow up limit is a
gradient shrinking soliton.

\end{abstract}
\maketitle
\tableofcontents
\section{Statement of result}

In the paper \cite{Ha:1} Hamilton introduced the concept of nonsingular solutions to a normalized Ricci flow and obtained important classification results on their infinite time limits in the 3 dimensional compact case. Let us recall that a suitably normalized Ricci flow is nonsingular if the curvature tensor is uniformly bounded for all time.
Around 2004 and 2006, with the help of Perelman's $W$ entropy \cite{P:1}, Sesum \cite{Se:1} and Fang, Zhang and Zhang \cite{FZZ:1} respectively  proved that noncollapsed, nonsingular solutions of certain normalized Ricci flows on compact manifolds in all dimensions sequentially converge to a gradient Ricci soliton as times goes to infinity. In a follow up paper \cite{FZZ:2}, the authors extended some of the results in \cite{FZZ:1} to the case of noncompact manifolds with finite volume. The proof of these results rely on the existence of a minimizer for the $W$ entropy which in turn depends crucially on the finiteness of the volume. In this note, we will prove a similar convergence result without any assumption on the volume or diameter of the underlying manifold.
A price to pay is that the convergence is a pointed one on selected points rather than all points.
Another feature is that we only need to assume the curvature tensor is bounded in a sequence of time intervals of fixed length.
As an application, for  a finite time singularity which is partially type I,  we prove that a blow up limit is a
gradient shrinking soliton.  Let us recall that a similar convergence to nonflat solitons for a type I singularity was proven around 10 years ago in \cite{EMT:1} and \cite{CZ:1} independently.  See also \cite{N:1}. In general such a statement, desirable for some researchers in the field,  is still lacking for a type II (non type I) singularity when the dimension is 4 or higher. Although recently in \cite{Ba:1}, a new type of convergence to singular solitons is presented for all finite time singularities. Also in 3 dimensional case, Perelman \cite{P:1} and \cite{P:2} already proved that a suitable blow up limit of a finite time singularity is a gradient Ricci soliton.

Before presenting the theorem,  let us fix the notations and definitions. In this note, $\M$ is a  Riemannian manifold of dimension $n$, which is either complete noncompact or compact without boundary.
Let $g=g(t)=g(\cdot, t)$, $t \in [0, T)$, $T>0$,  be  a  family of smooth metrics  on $\M$ and $Ric=Ric(t)=Ric(\cdot, t)$ be the  Ricci curvature corresponding to $g(t)$.  If the following equation holds:
\be
\lab{rf}
\partial_t g_{ij} = - 2 R_{ij},
\ee $t \in [0, T)$, then we say $(\M, g(t))$ is a  (un-normalized) Ricci flow on $[0, T)$;
if the equation
\be
\lab{nrf}
\partial_t g_{ij} = - 2 R_{ij} + g_{ij},
\ee holds for $t \in [0, T)$, then we say $(\M, g(t))$ is a  normalized Ricci flow on $[0, T)$.
We mentioned that if $\M$ is compact, the normalization constant in front of $g_{ij}$ is usually taken as
$\frac{2}{n} \oint R dg(t)$, where $R$ is the scalar curvature with respect to $g(t)$. Here, since the volume can be infinite, $1$ is chosen as the normalization constant. Under the  metric $g(t)$, we will  use $B(x, r, t)$ to denote the geodesic ball centered at $x$ with radius $r$ and $|B(x, r, t)|$ its volume; $Rm(\cdot, t)$ is the curvature tensor whose maximum norm at $x$ is $|Rm(x, t)|_{g(t)}$. Sometime, if there is no confusion, the time variable $t$ or $g(t)$ will be dropped from these quantities. In this note, we will always assume that $|Rm(\cdot, t)|_{g(t)}$ is bounded for each time $t$. But the bounds are not required to be uniform in time.

\begin{theorem}
\lab{thm1.1}
(a). Let $(\M, g(t))$ be a  normalized Ricci flow \eqref{nrf} on $[0, \infty)$, where $\M$  is either complete noncompact or compact without boundary.  Let $\{ t_j \}$ be  a  sequence of time going to infinity as $j \to \infty$. Suppose, there exist positive constants $\alpha, \beta$ and $\delta$ such that the following assumptions hold on the time interval $[t_j-\delta, t_j]$.
1.
$
| Rm | \le \alpha;
$
2. for all $x
\in \M$,
$
|B(x, 1, t)|_{g(t)} \ge \beta.
$

Then, there exists a subsequence of time, still denoted by $\{t_j\}$ and points $x_j \in \M$, such that
the Ricci flows $(\M, g(t_j-\delta +s), x_j)$, $s \in [0, \delta]$, converges, in $C^\infty_{loc}$ topology, to a gradient Ricci soliton.

(b).  Let $T>0$ be a finite number and $(\M, g(t))$ be a (un-normalized) Ricci flow \eqref{rf} on $[0, T)$, which is $\kappa$ non-collapsed and  which forms a singularity at $T$.
Let $\{ t_j \}$ be  a  sequence of time going up to $T$ as $j \to \infty$. Suppose, there exist positive constants $\alpha$ and $\delta$ such that the following assumption holds on the time interval $[t_j-\delta (T -t_j), t_j]$.
\[
| Rm | \le \alpha/(T-t_j).
\]Then, there exists a subsequence of time, still denoted by $\{t_j\}$ and points $x_j \in \M$, such that
the scaled Ricci flows $(\M, (T-t_j)^{-1} g(t_j-\delta (T-t_j) + (T-t_j) s), x_j)$, $s \in [0, \delta]$, converges, in $C^\infty_{loc}$ topology, to a gradient Ricci soliton.

\end{theorem}

Let us describe the idea of the proof.  As mentioned, the previous proof relies on the fact that the volume of the manifold is fixed along the  normalized Ricci flow and on the existence of
a minimizer for the $W$ entropy on compact manifolds.  Such existence
was proven in \cite{R:1}.
  The situation on noncompact manifolds is different.
  See \cite{Z12:1} for existence and nonexistence results of minimizers of the W entropy on noncompact manifolds.  Along the normalized Ricci flow \eqref{nrf}, the volume or diameter of the manifold may grow to infinity even in the compact setting. So even if a minimizer exists at each time $t$, it may decay to $0$  as $t \to \infty$.  To get around, we rely on an earlier observation that for the $W$ entropy with a fixed parameter, the maximum value of a minimizer on any compact domain  does not decay to $0$ as the domains fill up the whole manifold.  Then we show that for suitably chosen exhausting domains,  in a ball of fixed radius around a maximum point, the minimizer is bounded below by a fixed positive constant. Then we can apply Perelman's monotonicity formula for the $W$ entropy to reach the conclusion after taking limits. The lower bound is here to prevent a trivial $0=0$ limit from happening.

\section{proof}

\noindent Proof of Theorem \ref{thm1.1}.

(a). We start by recalling the entropies for the Ricci flow and their monotonicity property, introduced by Perelman \cite{P:1}.

Let $v \in W^{1, 2}({\M}, g(t)) \equiv \{v  \, | \,  \Vert v \Vert_{L^2(\M, g(t))} + \Vert \nabla v \Vert_{L^2(\M, g(t))} < \infty \}$ and $u=v^2$. Then
\be
\lab{defNv0}
N=N(v)(t)=\int_\M v^2 \ln v^2 dg = \int_\M u \ln u dg(t)
\ee is called Boltzmann's entropy; and
\be
\lab{defFv0}
F=F(v)(t) = \int_\M ( 4 |\nabla v |^2  + R v^2 )  dg = \int_\M \left( \frac{|\nabla u |^2}{u} + R u \right)dg(t)
\ee is Perelman's F entropy. Here $R$ is the scalar curvature of $(\M, g(t))$. Perelman's W entropy with parameter $\tau>0$ is
\be
\lab{defWv0}
\al
W&=W(t)= W(v, g(t), \tau) =\tau F(v) -N(v) -\left( \frac{n}{2} (\ln 4 \pi \tau) + n\right) \int_\M v^2 dg(t)\\
& =
 \int_\M \left[ \tau \left(\frac{|\nabla u |^2}{u} + R u \right)  - u \ln u  - \left(\frac{n}{2} (\ln 4 \pi \tau) + n \right) u \right] dg(t).
\eal
\ee

For the normalized Ricci flow \eqref{nrf}, the $W$ entropy with parameter $\tau=1$, which will be denoted by $W_1$ hereafter, plays a special role due to the following property.
Let $u=u(x, t)$ be a positive solution to the normalized conjugate heat equation
\be
\lab{nche}
\Delta u - (R - \frac{n}{2}) u + \p_t u =0, \qquad \M \times (0, T),
\ee such that $\int_\M u(x, t) dg(t) =1$. Then under suitable condition of $u$ and $\M$ at infinity, we have, for $v =u^2(\cdot, t)$, that
\be
\lab{dwdt}
\frac{d}{d t} W_1(t) \equiv \frac{d}{d t} W(v, g(t), 1) = 2 \int_\M \left| Ric + Hess \, \ln u - \frac{1}{2} g \right|^2(\cdot, t) \,  u(\cdot, t) d g(t).
\ee Here $Ric$ and $Hess $ are with respect to $g(t)$. If $\M$ is compact, the above formula is an immediate consequence of Perelman's monotonicity formula for the $W$ entropy of the un-normalized Ricci flow \cite{P:1}, after a change of variable. In the noncompact case, it is still valid if $(\M, g(t))$ has bounded curvature tensor and $u(x, t)$ satisfies
\be
\lab{vexpdecay}
u(x, t) \le C_1 e^{- C_2 d^2(x, x_0)}, \qquad x \in \M,
\ee where $C_1$ and $C_2$ are positive constants with $C_1$ depending on $|B(x_0, 1)|$. See Corollary 4.1 in \cite{Z12:1} e.g.
Note that the conditions can be improved to bounded Ricci curvature or even less, but we will not need the optimal one this time.

Given $t>0$,
let
\be
\lab{laminf}
\lambda(t)  \equiv \inf_{ } \{ W(v, g(t), 1)  \, | \, v \in C^\infty_0(\M), \,  \Vert v \Vert_{L^2(\M, g(t))} =1 \}
\ee be the infimum of the $W_1$ entropy with respect to $g(t)$. If $\M$ is compact, it is known for a long time from \cite{R:1} that this infimum is reached by some smooth function $v$. However, if $\M$
is noncompact, the infimum may not exist \cite{Z12:1}.  It is also clear from \eqref{dwdt} that
$\lambda(t)$ is nondecreasing in time even in the noncompact case, under the condition stated above. Indeed, fixing $t_0>0$, by definition of $\lambda(t_0)$, for any $\epsilon>0$, there exists, $v \in C^\infty_0(\M)$ with $ \Vert v \Vert_{L^2(\M, g(t_0))} =1$ such that
\be
\lab{wew}
W(v, g(t_0), 1) - \epsilon \le \lambda(t_0) \le W(v, g(t_0), 1).
\ee
Then we solve the conjugate heat equation with final value at time $t_0$.
\be
\lab{che2}
\al
\begin{cases}
&\Delta u - (R - \frac{n}{2}) u + \p_t u =0 , \quad \text{on} \quad \M \times [t_0-\delta, t_0];\\
&u(x, t_0) = v^2(x).
\end{cases}
\eal
\ee  Here $\delta$ is any small positive number. Then \eqref{dwdt} and \eqref{wew} infer that
\[
\lambda(t_0 - \delta) \le W(u^2(\cdot, t_0-\delta), g(t_0 - \delta), 1) \le W(v,  g(t_0), 1) \le
\lambda(t_0) + \epsilon.
\]Letting $\epsilon \to 0$, we see that $\lambda(t)$ is nondecreasing in time since $t_0>0$ is arbitrary and $\delta>0$ can be arbitrarily small.

For a domain $D$ in $\M$, let
\be
\lab{laminfd}
\lambda_D(t)  \equiv \inf_{ } \{ W(v, g(t), 1)  \, | \, v \in C^\infty_0(D), \,  \Vert v \Vert_{L^2(D, g(t))} =1 \}
\ee be the infimum of the $W_1$ entropy with respect to $g(t)$ on $D$.

Let $t_j$ be a sequence of time going to infinity as $j \to \infty$.  For each $j=1, 2, ...$, by definition of $\lambda(t_j)$, there exists a compact domain $D_j \subset \M$ such that
\be
\lab{lam1j}
\lambda_{D_j}(t_j) - \frac{1}{j} \le \lambda(t_j) \le \lambda_{D_j}(t_j).
\ee

 Since $D_j$ is compact, from \cite{R:1}, there exists a minimizer $0 \le v_j \in W^{1,2}_0(D_j) \cap C^\infty_0(D_j)$ so that
\[
\Vert v \Vert_{L^2(D_j, g(t_j))} =1, \qquad
\lambda_{D_j}(t_j) = W( v_j, g(t_j), 1)
\] and $v_j$ satisfies
\be
\lab{maineqvj}
 4 \Delta v_j - (R - \frac{n}{2}) v_j
+  2  v_j \ln v_j +   \lambda_{D_j}(t_j)  v_j  + c_n v_j= 0, \quad \text{on}  \quad D_j,
\ee where $c_n = \frac{n}{2} (\ln 4 \pi) + n $.  Here $\Delta$ and $R$ are with respect to the metric $g(t_j)$.

Since $v_j$ is $0$ on the boundary of $D_j$, a maximum point $x_j$ of $v_j$ must occur in the interior.
Hence
\[
\Delta v_j(x_j) \le 0.
\]This and \eqref{maineqvj} imply that
\be
\lab{vjc*}
v_j(x_j) = \max v_j \ge \exp \left([(\min R(\cdot, t_j)-(n/2))-\lambda_{D_j}(t_j) -c_n]/2 \right)
\ge c^*>0.
\ee Here we have used the fact that $\lambda_{D_j}(t_j)$ is bounded from above by a constant under assumptions 1 and 2. This can be checked by taking a test function in a unit ball contained in $D_j$.
The point is that the maximum value of $v_j$ is bounded away from $0$ for all $j=1, 2, ...$
Also, since for any $t>0$, there exists a $j$ such that $t \le t_j$, we also obtain
\be
\lab{lamjie}
\lam(t) \le \lam(t_j) \le \lam_{D_j}(t_j) \le C.
\ee

Next we will prove that the domains $D_j$ can be chosen so that for a fixed number $r^*>0$ and all $j=1, 2, ...$,
\be
\lab{vjc*/2}
v_j(x) \ge c^*/2, \quad x \in B(x_j, r^*, t_j) \subset D_j,
\ee that is $v_j$ is bounded away from $0$ in the ball of fixed radius.

For this purpose, let us recall the existence of a smooth, distance like function. Under our assumption of bounded curvature at time $t_j$, it is well known that there is
 a smooth function $L=L(x)$ on $\M$, which satisfies
\[
\al
&|\nabla L(x) | \le C_1, \qquad |\nabla^2 L(x) | \le C_1, \qquad x \in \M,\\
& d(x, 0) +1 \le L(x) \le  d(x, 0) + C_1.
 \eal
\] Here $C_1$ is a constant depending only on the curvature bound and dimension. Here, $\nabla, \nabla^2 $ and $d(x, 0)$ are with respect to $g(t_j)$.
For a proof, see Proposition 26.46 in
\cite{CCGGIIKLLN3:1} or \cite{Ta:1}.  We hereby choose
\[
D_j \equiv \{  x \in \M  \, | \, L(x) < l_j \}
\] where $l_j>0$ is a sufficiently large number so that \eqref{lam1j} holds and $l_j$ is not a critical value of $L$.  By Sard's theorem, such a $l_j$ exists. Thus $\partial D_j$ is a $C^2$ boundary which can be expressed by a uniform $C^2$ function locally in geodesic balls whose centers are at $\partial D_j$ and
  whose radius is $1/2$ of the injectivity radius of $\M$. The later is bounded from below uniformly by a positive constant by our assumption 1 and 2.  See  \cite{CGT:1}  and \cite{CLY:1} e.g.

Now we argue that $v_j$ is uniformly bounded in $C^\alpha$ norm, i.e., there exists a
 positive constant $C$ such that
 \be
 \lab{vkwinfty}
 \Vert v_j \Vert_{C^{\alpha}(D_j, g(t_j))} \le C.
 \ee A proof, similar to those in Section 2 of \cite{Z12:1}, goes as follows.
 We extend $v_j$ to a function on the whole manifold $\M$ by setting $v_j(x)=0$ when $x
 \in \M - D_j$. The extended function is still denoted by $v_j$. Then $v_j \in W^{1, 2}(\M)$, and $v_j$ satisfies the following inequality in the weak sense
 \[
 4 \Delta v_j - (R-\frac{n}{2}) v_j + 2 v_j \ln v_j + \lambda_{D_j}(t_j) v_j \ge 0, \qquad \text{in} \qquad \M.
 \]i.e. the extended $v_j$ is a subsolution.  Since $\lambda_{D_j}(t_j)$ and $R$ are uniformly bounded,  a  mean value inequality holds on unit balls, as given in Lemma 2.1 of \cite{Z12:1}. Namely
 \be
 \lab{mvi}
 \sup_{B(x, 1)} v^2_j \le C \int_{B(x, 2)}  v^2_j dg(t_j), \qquad x \in M.
 \ee Here $C$ depends only on the dimension, upper bound for $\lambda_{D_j}(t_j)$, the lower bound of the Ricci curvature and the lower bound of the volume of the unit balls.

   Using that $\Vert v_j \Vert_{L^2(M, g(t_j)}=1$, we see the norm
 $\Vert v_k \Vert_{L^\infty(\M)}$ is uniformly bounded. Hence the original $v_j$ in $D_j$ is actually a bounded solution
 to the Poisson equation
 \[
 \begin{cases}
 \Delta v_j(x) = f_j(x), \qquad x \in D_j \\
 v_j(x)=0, \quad x \in \partial D_j
 \end{cases}
 \] with $\Vert f_j \Vert_{L^\infty(D_j)} \le C$. Let $r_0$ be $1/2$ of the injectivity radius of $(\M, g(t_j))$. Let
 \[
 D_{j, r_0} \equiv \{ x \in D_j \, | \, d(x, \partial D_j)>r_0/4 \}.
 \]Then
the standard  interior regularity theory of elliptic equations shows $\Vert v_j \Vert_{C^{\alpha}(D_{j, r_0}, g(t_j))} \le C.$ If $x \in D_j/D_{j, r_0}$, then
\[
x \in B(y, r_0/2)
\]for some $y \in \partial D_j$. Since, by our construction, $\partial D_j \cap B(y, r_0)$ is given by a $C^2$ function with uniformly bounded $C^2$ norm, the standard boundary regularity theory of elliptic equations  infers that $\Vert v_j \Vert_{C^{\alpha}(B(y, r_0/2), g(t_j))} \le C.$  Since all points in $D_j$ are covered,  (\ref{vkwinfty}) is true. From \eqref{vjc*} and \eqref{vkwinfty}, we obtain \eqref{vjc*/2} for some fixed
$r^*>0$.

Next we solve the conjugate heat equation with final value at time $t_j$:
\be
\lab{chej}
\al
\begin{cases}
&\Delta u_j - (R - \frac{n}{2}) u_j + \p_t u_j =0 , \quad \text{on} \quad \M \times [t_j-\delta, t_j];\\
&u_j(x, t_j) = v^2_j(x).
\end{cases}
\eal
\ee  Here $\delta$ is any small positive number.  Fixing $y \in \M$  and $t_j$, let $G=G(x, t; x_j, t_j)$ be the heat kernel of the equation in\eqref{chej}.  Then for $x \in \M$ and $t \in [t_j -\delta, t_j]$,
\[
u_j(x, t) = \int_\M G(x, t; y, t_j) v^2_j(y) dg(t_j).
\]Under our assumption, it is well known that the following upper and lower bound for $G$ holds on the time interval $[t_j-\delta, t_j]$. See Chapter 26 of \cite{CCGGIIKLLN3:1} e.g.
\[
\frac{1}{C_1 (t_j-t)^{n/2}} e^{- d^2(x, y, t_j)/[c_2(t_j-t)]} \le G(x, t; y, t_j) \le \frac{C_1}{(t_j-t)^{n/2}} e^{-c_2 d^2(x, y, t_j)/(t_j-t)}.
\] We comment that these are short time generic bounds which hold much less curvature assumption than that of bounded curvature. Since $v_j$ are compactly supported and bounded, the upper bound for $G$ implies quadratic
exponential decay for $u_j$.  On the other hand, the lower bound for $G$ and \eqref{vjc*} tells us, for $t \in [t_j-\delta, t_j]$,
\be
\lab{uj>fjl}
\al
u_j(x, t) &\ge (c^*/2)^2 \int_{B(x_j, r^*)} \frac{1}{C_1 (t_j-t)^{n/2}} e^{- d^2(x, y, t_j)/[c_2(t_j-t)]} dg(t_j) \\
&\ge c_3 e^{-c_4 [1+d^2(x, x_j, t_j)]/(t_j-t)}.
\eal
\ee Here $c_3, c_4$ are positive constants independent of $j$.

As explained earlier, the fast decay of $u_j$ allows us to justify  \eqref{dwdt} which gives us
\[
\al
 W(v_j, g(t_j), 1) -& W(u^2(\cdot, t_j - \delta),  g(t_j), 1) \\
  &= 2 \int^{t_j}_{t_j - \delta} \int_\M \left| Ric + Hess \, \ln u_j - \frac{1}{2} g \right|^2(\cdot, t) \,  u_j(\cdot, t) d g(t) dt.
 \eal
\]From this and \eqref{lam1j} we see that
\be
\lab{lamjjdel}
 \frac{1}{j} + \lambda(t_j) - \lam(t_j-\delta)
   \ge 2 \int^{t_j}_{t_j - \delta} \int_\M \left| Ric + Hess \, \ln u_j - \frac{1}{2} g \right|^2(\cdot, t) \,  u_j(\cdot, t) d g(t) dt.
\ee  According to Hamilton's compactness theorem \cite{Ha:2}, the Ricci flows
\[
(\M, g(t_j-\delta+s), x_j), \qquad s \in [0, \delta],
\]converge sub-sequentially in $C^\infty_{loc}$ topology to a smooth Ricci flow $(\M_\infty, g_\infty(s), x_\infty)$. The boundedness of $u_j$ and \eqref{uj>fjl} implies that $u_j(\cdot, \cdot+t_j-\delta)$ converges in $C^\alpha_{loc}$ topology to a  function $u_\infty$ which is positive on $\M_\infty \times (0, \delta]$. We mention that the uniform boundedness of $u_j$ is a simple consequence of the boundedness of $v_j$ and the heat kernel upper bound.  Compactness in $C^\alpha_{loc}$ topology comes from the H\"older continuity of $v_j$ and gradient bound of heat kernel.

As shown earlier in \eqref{lamjie}, $\lam=\lam(t)$ is a bounded, nondecreasing function. Hence
\[
\lim_{j \to \infty} (\lam(t_j)-\lam(t_j-\delta)) =0.
\] This, \eqref{lamjjdel} and Fatou's lemma imply that
\[
0
=\int^{\delta}_0 \int_{\M_\infty} \left| Ric + Hess \, \ln u_\infty - \frac{1}{2} g_\infty \right|^2(\cdot, s) \,  u_\infty(\cdot, s) d g_\infty(s) ds.
\]Hence $(\M_\infty, g_\infty(s), x_\infty)$ is a gradient shrinking soliton. This proves part (a) of the theorem.

Next we prove part (b).  Without loss of generality, we take $T=1$.  For the un-normalized Ricci flow $g_{ij}(t)$ given in part (b), we carry out the scaling
\[
s = - \ln (1-t), \qquad \tilde{g}_{ij}(s) =\frac{1}{1-t} g_{ij}(t).
\]Then we can apply part (a) on the normalized Ricci flow $\tilde{g}_{ij}(s)$ to complete the proof.
Note that the $\kappa$ non-collapsing condition for $g_{ij}(t)$ and the curvature bound infer
the volume non-collapsing assumption in part (a). \qed

{\bf Acknowledgment.}
Support by the Simons Foundation  grant 710364 is gratefully acknowledged.
 The author also wishes to thank  Prof. Zhu Meng for helpful conversations.

\bigskip

\noindent e-mail:  qizhang@math.ucr.edu

\enddocument